\documentclass[12pt]{article}

\usepackage{cite}
\usepackage{amsfonts, amssymb, amsmath, amsthm, epsf}

\theoremstyle{plain}
\newtheorem{theorem}{Theorem}

\newtheorem{condition}{Condition}

\theoremstyle{definition}

\newtheorem{remark}{Remark}

\newcommand{\supp}{{\rm supp\,}}
\newcommand{\spt}{{\rm spt\,}}
\newcommand{\Dom}{{\rm D}}

\newcommand{\dist}{{\rm dist}}
\renewcommand{\phi}{{\varphi}}
\newcommand{\n}{{|\!|\!|}}

\newcommand{\cK}{{\mathcal K}}
\newcommand{\cM}{{\mathcal M}}
\newcommand{\cN}{{\mathcal N}}

\newcommand{\cO}{{\mathcal O}}

\newcommand{\cR}{{\mathcal R}}

\newcommand{\cF}{{\mathcal F}}

\newcommand{\pG}{{\partial G}}

\newcommand{\oG}{{\overline G}}

\newcommand{\bP}{{\mathbf P}}
\newcommand{\bB}{{\mathbf B}}

\newcommand{\bbR}{{\mathbb R}}

\begin{document}

%Series: Mahtematics\hfill UDC 517.9

%\baselineskip=21pt

\begin{center}
{\bf\Large Unbounded perturbations of two-dimensional diffusion
processes with nonlocal boundary conditions}

\smallskip

Pavel Gurevich\footnote{This research was supported by Russian
Foundation for Basic Research (project No.~07-01-00268) and the
Alexander von Humboldt Foundation.}
\end{center}

\bigskip

\abstract{The existence of Feller semigroups arising in the theory
of multidimensional diffusion processes is studied. Unbounded
perturbations of elliptic operators (in particular,
integro-differential operators) are considered in plane bounded
regions. Their domain of definition is given by a nonlocal
boundary condition involving the integral over the closure of the
region with respect to a nonnegative Borel measure. The support of
the measure may intersect with the boundary, and the
 measure need not be small.  We
formulate sufficient conditions on unbounded perturbations of
elliptic operator and on the Borel measure in the nonlocal
boundary condition which guarantee that the corresponding nonlocal
operator is a generator of a Feller semigroup.}

\bigskip

{\bf 1.}  A general form of a generator of a strongly continuous contractive nonnegative
semigroup (Feller semigroup) of operators acting between the spaces of continuous
functions on an interval was investigated in~\cite{Feller2}. In the multidimensional
case, it was proved that the generator of a Feller semigroup is an elliptic differential
operator (possibly with degeneration) whose domain of definition consists of continuous
functions satisfying nonlocal conditions which involve an integral over the closure of
the region with respect to a nonnegative Borel measure~\cite{Ventsel}. The inverse
problem remains open: given an elliptic integro-differential operator whose domain of
definition is described by nonlocal boundary conditions, whether or not the closure of
this operator is a generator of a Feller semigroup. In~\cite{SU, BCP, Taira3,  Ishikawa,
GalSkJDE}, the authors considered the   {\it transversal\/} case, where the order of
nonlocal terms is less than the order of local terms. The {\it nontransversal\/} case,
where these orders coincide, was studied in~\cite{SkubDAN89,GalSkJDE} (see also the
bibliography in~\cite{GalSkJDE}).

It was assumed in~\cite{GalSkJDE} that the coefficients at the nonlocal terms are less
than one. In this paper, we  consider nontransversal nonlocal conditions on the boundary
of a plane domain and investigate the ``limit case,'' where the coefficients at nonlocal
terms may equal one at some points (the Borel measure is assumed to be atomic at these
points). We formulate sufficient conditions on unbounded perturbations of elliptic
operator and on the Borel measure in nonlocal conditions which guarantee that the
corresponding nonlocal operator is a generator of a Feller semigroup.

\smallskip

Let $G\subset\bbR^2$ be a bounded domain with boundary $\pG$.
Consider a set ${\cK}\subset\partial G$ consisting
of finitely many points. Let $\partial G\setminus{\mathcal
K}=\bigcup_{i=1}^{N}\Gamma_i$, where $\Gamma_i$ are open (in the
topology of $\partial G$) $C^\infty$ curves. Assume that the
domain $G$ is a plane angle in some neighborhood of each point
$g\in{\mathcal K}$.

For any closed sets $Q\subset\oG$ and $K\subset\oG$ such that $Q\cap K\ne\varnothing$, we
introduce the space $ C_K(Q)=\{u\in C(Q): u(y)=0,\ y\in Q\cap K\} $ with the
maximum-norm. Introduce the space $H_a^k(G)$ as a completion of the set of infinitely
differentiable functions vanishing near $\mathcal K$ with respect to the norm
$$
 \|u\|_{H_a^k(G)}=\Bigg(
    \sum_{|\alpha|\le k}\int_G \rho^{2(a+|\alpha|-k)} |D^\alpha u(y)|^2 dy
                                       \Bigg)^{1/2},
$$
where $a\in \mathbb R$, $k\ge 0$ is an integer, and
$\rho=\rho(y)=\dist(y,\mathcal K)$. For an integer $k\ge1$, denote
by $H_a^{k-1/2}(\Gamma)$  the space of traces of the functions
from $H_a^k(G)$ on a smooth curve $\Gamma\subset\overline{G}$. We
also introduce the norms in weighted spaces depending on the
parameter $q>0$. Set $
\|v\|_{H_{a}^0(\Gamma_i)}=\big(\,\int_{\Gamma_i}\rho^{2a}|v(y)|^2d\Gamma\big)^2,
$
$$
\n u \n_{H_a^k(G)}=
\left(\|u\|_{H_{a}^k(G)}^2+q^k\|u\|^2_{H_a^0(G)}\right)^{1/2},\qquad
k\ge0,
$$
$$
\n v \n_{H_a^{k-1/2}(\Gamma_i)}=
\left(\|v\|_{H_a^{k-1/2}(\Gamma_i)}^2+q^{k-1/2}\|
v\|_{H_{a}^0(\Gamma_i)}^2\right)^{1/2},\qquad k\ge1.
$$

{\bf 2.} Consider the differential operator $
P_0u=\sum_{i,j=1}^{2}p_{ij}(y)u_{y_iy_j}(y)+
\sum_{i=1}^2p_i(y)u_{y_i}(y)+p_0(y)u(y), $ where $p_{ij},p_i\in
C^\infty(\bbR^2)$ are real-valued functions, $p_0\ge0$, and
$p_{ij}=p_{ji}$, $i,j=1,2$. We  assume that there is a constant
$c>0$ such that $\sum_{i,j=1}^{2}p_{ij}(y)\xi_i\xi_j\ge c|\xi|^2$
for $y\in\overline{G}$ and $\xi=(\xi_1,\xi_2)\in\bbR^2.$

Introduce  operators corresponding to nonlocal terms supported near the set $\mathcal K$.
For any  set $\mathcal M$, we denote its $\varepsilon$-neighborhood by $\mathcal
O_{\varepsilon}(\mathcal M)$. Let $\Omega_{is}$ ($i=1, \dots, N;$ $s=1, \dots, S_i$) be
$C^\infty$ diffeomorphisms taking some neighborhood ${\mathcal O}_i$ of the curve
$\overline{\Gamma_i\cap\mathcal O_{{\varepsilon}}(\mathcal K)}$ to the set
$\Omega_{is}({\mathcal O}_i)$ in such a way that $\Omega_{is}(\Gamma_i\cap\mathcal
O_{{\varepsilon}}(\mathcal K))\subset G$ and $ \Omega_{is}(g)\in\mathcal K$ for $
g\in\overline{\Gamma_i}\cap\mathcal K. $ Thus, the transformations $\Omega_{is}$ take the
curves $\Gamma_i\cap\mathcal O_{{\varepsilon}}(\mathcal K)$  inside the domain $G$ and
the set of their end points $\overline{\Gamma_i}\cap\mathcal K$ to itself.

Let us specify the structure of the transformations $\Omega_{is}$. Denote by $\Omega_{is}^{+1}$ the
transformation $\Omega_{is}:{\mathcal O}_i\to\Omega_{is}({\mathcal
O}_i)$ and by $\Omega_{is}^{-1}:\Omega_{is}({\mathcal
O}_i)\to{\mathcal O}_i$ the inverse transformation. The set of
points
$\Omega_{i_qs_q}^{\pm1}(\dots\Omega_{i_1s_1}^{\pm1}(g))\in{\mathcal
K}$ ($1\le s_j\le S_{i_j},\ j=1, \dots, q$) is said to be an {\em
orbit} of the point $g\in{\mathcal K}$. In other words, the orbit
of a point $g\in\cK$ is formed by the points (of the set $\mathcal
K$) that can be obtained by consecutively applying the
transformations $\Omega_{i_js_j}^{\pm1}$ to the point $g$. The set
$\mathcal K$ consists of finitely many disjoint orbits, which we
denote by $\mathcal K_\nu$, $\nu=1,\dots,N_0$.

Take a sufficiently small number $\varepsilon>0$ such that there
exist neighborhoods $\mathcal O_{\varepsilon_1}(g_j)$, $ \mathcal
O_{\varepsilon_1}(g_j)\supset\mathcal O_{\varepsilon}(g_j) $,
satisfying the following conditions: 1. the domain $G$ is a plane
angle in the neighborhood $\mathcal O_{\varepsilon_1}(g_j)$;  2.
$\overline{\mathcal O_{\varepsilon_1}(g)}\cap\overline{\mathcal
O_{\varepsilon_1}(h)}=\varnothing$ for any $g,h\in\mathcal K$,
$g\ne h$; 3. if $g_j\in\overline{\Gamma_i}$ and
$\Omega_{is}(g_j)=g_k,$ then ${\mathcal
O}_{\varepsilon}(g_j)\subset\mathcal
 O_i$ and
 $\Omega_{is}\big({\mathcal
O}_{\varepsilon}(g_j)\big)\subset{\mathcal
O}_{\varepsilon_1}(g_k).$

For each point $g_j\in\overline{\Gamma_i}\cap\mathcal K_\nu$, we
fix a linear transformation $Y_j: y\mapsto y'(g_j)$ (the
composition of the shift by the vector $-\overrightarrow{Og_j}$
and rotation) mapping the point $g_j$ to the origin in such a way
that $ Y_j({\mathcal O}_{\varepsilon_1}(g_j))={\mathcal
O}_{\varepsilon_1}(0),\ Y_j(G\cap{\mathcal
O}_{\varepsilon_1}(g_j))=K_j\cap{\mathcal O}_{\varepsilon_1}(0), $
$ Y_j(\Gamma_i\cap{\mathcal
O}_{\varepsilon_1}(g_j))=\gamma_{j\sigma}\cap{\mathcal
O}_{\varepsilon_1}(0)\ (\sigma=1\ \text{or}\ 2), $ where $
 K_j$ is a plane angle of nonzero opening and $\gamma_{j\sigma}$ its sides.

\begin{condition}\label{condK1}
Let $g_j\in\overline{\Gamma_i}\cap\mathcal K_\nu$ and
$\Omega_{is}(g_j)=g_k\in\mathcal K_\nu;$ then the transformation $
Y_k\circ\Omega_{is}\circ Y_j^{-1}:{\mathcal
O}_{\varepsilon}(0)\to{\mathcal O}_{\varepsilon_1}(0) $ is the
composition of rotation and homothety.
\end{condition}

Introduce the nonlocal operators: $
 \mathbf B_{i}u=\sum_{s=1}^{S_i}
b_{is}(y) u(\Omega_{is}(y))$, $
   y\in\Gamma_i\cap\mathcal O_{\varepsilon}(\mathcal K),$ and $
\mathbf B_{i}u=0,$  $ y\in\Gamma_i\setminus\mathcal
O_{\varepsilon}(\mathcal K),$ where $b_{is}\in C^\infty(\mathbb
R^2)$ are real-valued functions, $\supp b_{is}\subset\mathcal
O_{{\varepsilon}}(\mathcal K)$.

\begin{condition}\label{cond1.2}
\begin{enumerate}
\item
$ b_{is}(y)\ge0,\quad \sum_{s=1}^{S_i} b_{is}(y)\le 1,\quad
 y\in\overline{\Gamma_i};
$
\item
$
 \sum_{s=1}^{S_i}
b_{is}(g)+\sum_{s=1}^{S_j} b_{js}(g)<2,\
g\in\overline{\Gamma_i}\cap\overline{\Gamma_j}\subset\cK,\quad\text{if}\
 i\ne j\ \text{and}\
 \overline{\Gamma_i}\cap\overline{\Gamma_j}\ne\varnothing.
$
\end{enumerate}
\end{condition}

\begin{theorem}\label{th} Let Conditions~$\ref{condK1}$ and $\ref{cond1.2}$ be fulfilled.
Then there is a number $\delta_0>0$ such that, for
$k=0,1,2,\dots$, $\delta\in(0,\delta_0)$,   $q> q_0(\delta)\ge 0$,
and $\psi_i\in C_\cK(\overline{\Gamma_i})\cap
H_{a}^{k+3/2}(\Gamma_i)$, $a=k+1-\delta$, the problem
$$
P_0u-q u =0, \  y\in G;\qquad u|_{\Gamma_i}-\bB_i u=\psi_i(y), \
y\in\Gamma_i,\ i=1,\dots,N,
$$
has a unique solution $u\in C^\infty(G)\cap C_\cK(\overline G)\cap H_{a}^{k+2}(G);$
moreover,
$$\|u\|_{C(\oG)}+\n u\n_{H_{a}^{k+2}(G)}\le c
\sum_{i=1}^N\big(\|\psi_i\|_{C(\overline{\Gamma_i})}+\n\psi_i\n_{H_{a}^{k+3/2}(\Gamma_i)}\big),$$
where $c>0$ does not depend  on  $q$, $\psi_i$, or $u$.
\end{theorem}

From now on, we fix $\delta\in(0,\delta_0)$, an integer $k\ge2$,
and $a=k+1-\delta$.

{\bf 3.} Consider a linear bounded operator $P_1: H_a^{k+2}(G)\to
H_{a-1}^k(G)$.

\begin{condition}\label{cond3.1}
\begin{enumerate}
\item
If a function $u\in H_a^{k+2}(G)$ achieves its positive maximum at
a point $y^0\in G$, then $P_1 u(y^0)\le0$.
\item
If $u\in C(\oG)\cap H_a^{k+2}(G)$, then the function $P_1u$ is
bounded on $G$.
\item
For all sufficiently small $\varrho>0$, we have $
P_1=P_{1\varrho}^1+P_{1\varrho}^2, $ where the operators
$P_{1\varrho}^1,P_{1\varrho}^2:H_a^{k+2}(G)\to H_{a-1}^k(G)$ are
such that
\begin{enumerate}
\item
$\| P_{1\varrho}^1u\|_{H_{a-1}^k(G)}\le c(\varrho)\|
u\|_{H_a^{k+2}(G)}$, $c(\varrho)>0$ does not depend on $u$ and
$c(\varrho)\to0$ as $\varrho\to 0$,
\item
the operator $P_{1\varrho}^2$ is compact.
\end{enumerate}
\end{enumerate}
\end{condition}
Note that $\Dom(P_1)\subset C^k(\oG\setminus\cK)\subset C^2(G)$
and $\cR(P_1)\subset C^{k-2}(\oG\setminus\cK)\subset C(G)$.
However, $\cR(P_1)\not\subset C(\oG)$ in general.

\begin{remark}
The prototype of the abstract operators $P_1$, $P_{1\varrho}^1$, and $P_{1\varrho}^2$ are
integral operators of the form  (cf.~\cite{GalSkJDE, Taira3,GarMenaldi})
\begin{equation*}
\begin{aligned}
 P_{1}u(y)&=\int_{F} [u(y+z(y,\eta))-u(y)-(\nabla
u(y),z(y,\eta))]m(y,\eta)\pi(d\eta),\\
 P_{1\varrho}^1u(y)&=\int_{Z\le\varrho} [u(y+z(y,\eta))-u(y)-(\nabla
u(y),z(y,\eta))]m(y,\eta)\pi(d\eta),\\
P_{1\varrho}^2u(y)&=\int_{Z>\varrho} [u(y+z(y,\eta))-u(y)-(\nabla
u(y),z(y,\eta))]m(y,\eta)\pi(d\eta),
\end{aligned}
\end{equation*}
where  $F$ is a space with a $\sigma$-algebra $\cF$ and a Borel measure $\pi$,
 $ y+  z(y,\eta)\in\oG$ and $|D^\alpha_y z(y,\eta)|\le Z(\eta)$ for  $y\in\oG,\ \eta\in F$, $|\alpha|\le k$, $Z(\eta)$
is a nonnegative $\pi$-measurable bounded function, and  $m(y,\eta)\ge0$ (some additional
restrictions on the functions $z(y,\eta)$, $Z(\eta)$, and $m(y,\eta)$ should also be
imposed).
\end{remark}

{\bf 4.} In this paper, we consider the nonlocal conditions in the
 {\it nontransversal} case (see, e.g.,~\cite{Taira3} for the probabilistic interpretation):
\begin{equation}\label{eq57}
u(y)-\int_\oG u(\eta)\mu_i(y,d\eta)=0,\ y\in\Gamma_i,\
i=1,\dots,N;\qquad u(y) =0,\  y\in\cK,
\end{equation}
where $\mu_i(y,\cdot)$ is a nonnegative Borel measure on $\oG$
such that $\mu_i(y,\oG)\le 1,\ y\in\Gamma_i$.

Introduce the measures $\delta_{is}(y,\cdot)$ as follows: $
\delta_{is}(y,Q)=b_{is}(y)\chi_Q(\Omega_{is}(y)),$
$y\in\Gamma_i\cap\cO_\varepsilon(\cK),$ and $\delta_{is}(y,Q)=0$, $
y\in\Gamma_i\setminus\cO_\varepsilon(\cK)$, where $Q\subset\oG$ is an arbitrary Borel set
and $\chi_Q(\cdot)$ a characteristic function of the set $Q$.

 We study those  measures $\mu_i(y,\cdot)$ which
can be represented in the form
$$
\mu_i(y,\cdot)=\sum_{s=1}^{S_i}\delta_{is}(y,\cdot)+\alpha_i(y,\cdot)+\beta_i(y,\cdot),\qquad
y\in\Gamma_i,
$$
where $\alpha_i(y,\cdot)$ and $\beta_i(y, \cdot)$ are nonnegative
Borel measures to be specified below.

Denote $ \spt\alpha_i(y,\cdot)=\oG\setminus\bigcup_{V\in T}\{V\in
T: \alpha_i(y,V\cap\oG)=0\} $ ($T$ is the set of all open sets in
$\bbR^2$).

Set $ \bB_{\alpha i}u(y)=\int_\oG u(\eta)\alpha_i(y, d\eta),\
\bB_{\beta i}u(y)=\int_\oG u(\eta)\beta_i(y, d\eta),\
y\in\Gamma_i. $

We assume that the measures $\alpha_i(y,\cdot)$ and $\beta_i(y,\cdot)$ satisfy the
following conditions (cf.~\cite{GalSkJDE}).

\begin{condition}\label{cond2.3}
There exist numbers $\varkappa_1>\varkappa_2>0$ and $\sigma>0$
such that
\begin{enumerate}
\item
$\spt\alpha_i(y,\cdot)\subset\oG\setminus\cO_{\varkappa_1}(\cK)$,
  $y\in\Gamma_i$, and $ \n
\bB_{\alpha i}u\n_{H_a^{k+3/2}(\Gamma_i)}\le c\n
u\n_{H_a^{k+2}(G\setminus\overline{\cO_{\varkappa_1}(\cK)})}, $
\item
$\spt\alpha_i(y,\cdot)\subset\overline{G_\sigma}$,
$y\in\Gamma_i\setminus\cO_{\varkappa_2}(\cK)$, and $ \n
\bB_{\alpha
i}u\n_{H_a^{k+3/2}(\Gamma_i\setminus\overline{\cO_{\varkappa_2}(\cK)})}\le
c\n u \n_{H_a^{k+2}(G_\sigma)}, $
\end{enumerate}
where $G_\sigma=\{y\in G:\dist(y,\pG)<\sigma\}.$
\end{condition}

Let $
 \cN=\bigcup_{i=1}^N\{y\in\Gamma_i: \mu_i(y,\oG)=0\}\cup\cK$ and $\cM=\pG\setminus \cN.$
Assume that $\cN$ and $\cM$ are Borel sets.

\begin{condition}\label{cond2.4}
 $\beta_i(y,\cM)<1$ for $y\in\Gamma_i\cap\cM$, $i=1,\dots, N$.
\end{condition}

\begin{condition}\label{cond2.5}
For any function $u\in C_\cN(\overline G)$, the functions
$\bB_{\alpha i}u$ and $\bB_{\beta i}u$ can be extended to
$\overline{\Gamma_i}$ in such a way that the extended functions
{\rm(}which we also denote by $\bB_{\alpha i}u$ and $\bB_{\beta
i}u$, respectively{\rm)} belong to $C_\cN(\overline{\Gamma_i})$.
\end{condition}

We represent the measures $\beta_i(y,\cdot)$ in the form $
\beta_i(y,\cdot)=\beta_i^1(y,\cdot)+\beta_i^2(y,\cdot), $ where
$\beta_i^1(y,\cdot)$ and $\beta_i^2(y,\cdot)$ are nonnegative
Borel measures. Let us specify them.  Set
$\cM_p=\cO_p(\overline{\cM})$, $p>0$. Consider a cut-off function
$\hat\zeta_p\in C^\infty(\bbR^2)$ such that $0\le\hat\zeta_p(y)\le
1$, $\hat\zeta_p(y)=1$ for $y\in\cM_{p/2}$, and $\hat\zeta_p(y)=0$
for $y\notin\cM_{p}$. Set $\tilde\zeta_p=1-\hat\zeta_p$. Introduce
the operators
$$
\hat\bB_{\beta i}^1
u(y)=\int_{\oG}\hat\zeta_p(\eta)u(\eta)\beta_i^1(y,d\eta),\qquad
 \tilde\bB_{\beta i}^1 u(y)=\int_{\oG}\tilde\zeta_p(\eta)
u(\eta)\beta_i^1(y,d\eta),$$ $$  \bB_{\beta i}^2
u(y)=\int_{\oG}u(\eta)\beta_i^2(y,d\eta).
$$

Consider the  Banach spaces with the norms depending on the parameter $q>0$:
\begin{enumerate}
\item
 $H_{\cN,a}^{k+2}(G)=C_\cN(\oG)\cap
H_a^{k+2}(G)$, $ \n u
\n_{H_{\cN,a}^{k+2}(G)}=\|u\|_{C_\cN(\oG)}+\n u \n_{H_a^{k+2}(G)},
$
\item
 $H_{\cN,a}^{k+3/2}(\Gamma_i)=C_\cN(\overline{\Gamma_i})\cap
H_a^{k+3/2}(\Gamma_i)$, $ \n v
\n_{H_{\cN,a}^{k+3/2}(\Gamma_i)}=\|v\|_{C_\cN(\overline{\Gamma_i})}+\n
v \n_{H_a^{k+3/2}(\Gamma_i)}.$
\end{enumerate}

\begin{condition}\label{cond3.3}
For each sufficiently small $p>0$ the following assertions are
true{\rm:}
\begin{enumerate}
\item
 $\n \hat\bB_{\beta i}^1\n_{H_{\cN,a}^{k+2}(G)\to
H_{\cN,a}^{k+3/2}(\Gamma_i)}\to 0,\  p\to0,$ uniformly with
respect to $q;$
\item
the norms $\n \tilde\bB_{\beta i}^1\n_{H_{\cN,a}^{k+2}(G)\to
H_{\cN,a}^{k+3/2}(\Gamma_i)}$ are bounded uniformly with respect
to $q$.
\end{enumerate}
\end{condition}

\begin{condition}\label{cond3.4}
The operators $\bB_{\beta i}^2:H_{\cN,a}^{k+2}(G)\to
H_{\cN,a}^{k+3/2}(\Gamma_i)$, $i=1,\dots,N$, are compact.
\end{condition}

Introduce the space  $ C_B(\oG)=\{u\in C(\oG): u\ \text{satisfy
nonlocal conditions \eqref{eq57}}\}. $ We consider the unbounded
operator $\bP : \Dom(\bP )\subset C_B(\overline G)\to
C_B(\overline G)$ given by
$$
\bP u=P_0u+P_1u,\qquad u\in \Dom(\bP )=\{u\in C_B(\oG)\cap
H_{a}^{k+2}(G): P_0u+P_1u\in C_B(\overline G)\}.
$$
Note that $\Dom(\bP)\subset C^2(G)\cap C_B(\oG)$ due to the relation $k\ge 2$ and the
Sobolev embedding theorem. The proof of the following main result is based on
Theorem~\ref{th},  the Hille--Iosida theorem, and the maximum principle.
\begin{theorem}\label{th3.1}
Let Conditions~$\ref{condK1}$--$\ref{cond3.4}$ hold. Then the
operator $\bP$ admits the closure
$\overline{\bP}:\Dom(\overline{\bP})\subset C_B(\oG)\to C_B(\oG)$
and the operator $\overline{\bP}$  is a generator of a Feller
semigroup.
\end{theorem}

{\bf 6.} Consider an example of nonlocal conditions in which the
assumptions of the paper hold. Let
$\pG=\Gamma_1\cup\Gamma_2\cup\cK$, where $\Gamma_1$ and $\Gamma_2$
are $C^\infty$ curves open and connected in the topology of $\pG$,
$\Gamma_1\cap\Gamma_2=\varnothing$, and
$\overline{\Gamma_1}\cap\overline{\Gamma_2}=\cK$; the set $\cK$
consists of two points $g_1$ and $g_2$.   Let $\Omega_j$,
$j=1,\dots,4$, be nondegenerate transformations of class $C^{k+2}$
defined on a neighborhood of $\overline{\Gamma_1}$ and satisfying
the following conditions (see Fig.~\ref{figEx4-4}):
\begin{figure}[ht]
{ \hfill\epsfxsize120mm\epsfbox{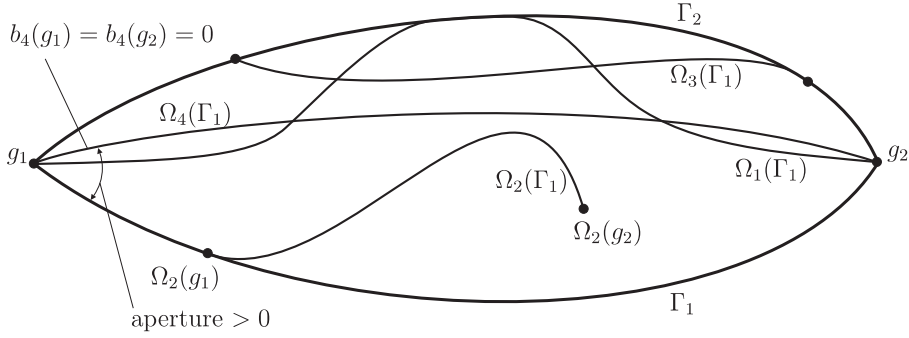}\hfill\ }
\caption{Nontransversal nonlocal conditions}
   \label{figEx4-4}
\end{figure}
\begin{enumerate}
\item
$\Omega_1(\cK)=\cK$,
$\Omega_1(\Gamma_1\cap\cO_\varepsilon(\cK))\subset G$,
$\Omega_1(\Gamma_1\setminus\cO_\varepsilon(\cK))\subset
G\cup\Gamma_2$, and $\Omega_1(y)$ is a composition of a shift of
the argument, rotation, and homothety  for
$y\in\overline{\Gamma_1}\cap\cO_\varepsilon(\cK)$;
\item
there exist numbers $\varkappa_1>\varkappa_2>0$ and $\sigma>0$
such that
$\Omega_2(\overline{\Gamma_1})\subset\oG\setminus\cO_{\varkappa_1}(\cK)$
and
$\Omega_2(\overline{\Gamma_1}\setminus\cO_{\varkappa_2}(\cK))\subset
\overline{G_\sigma}$; moreover, $\Omega_2(g_1)\in\Gamma_1$ and
$\Omega_2(g_2)\in G$;
\item
$\Omega_3(\overline{\Gamma_1})\subset G\cup{\Gamma_2}$ and
$\Omega_3(\cK)\subset\Gamma_2$;
\item
$\Omega_4(\overline{\Gamma_1})\subset G\cup\overline{\Gamma_2}$
and $\Omega_4(\cK)=\cK$; the angle between the rays tangent  to
$\Gamma_1$ and $\Omega_4(\Gamma_1)$ at the point $g_j$ is nonzero.
\end{enumerate}

Let  $b_j\in C^{k+2}(\overline{\Gamma_1})$ and $b_j\ge0$, $j=1,\dots,4$.

Consider the following nonlocal conditions:
$$
\begin{aligned}
u(y)-\sum_{j=1}^4 b_j(y)u(\Omega_j(y)) &=0, &
&y\in\Gamma_1,\\
u(y) &=0, & &y\in\overline{\Gamma_2}.
\end{aligned}
$$

We assume that
\begin{equation*}
\begin{gathered}
\sum_{j=1}^4 b_j(y)\le 1,\ y\in\Gamma_1;\\
 b_2(g_1)=0\
\text{or}\ \sum_{j=1}^4 b_j(\Omega_2(g_1))=0;\quad  b_2(g_2)=0;\quad b_4(g_j)=0.
\end{gathered}
\end{equation*}

Introduce a cut-off function $\zeta\in C^\infty(\bbR^2)$ supported
in $\cO_\varepsilon(\cK)$, equal to $1$ on
$\cO_{\varepsilon/2}(\cK)$, and such that $0\le\zeta(y)\le 1$,
 $y\in\bbR^2$. Let
$y\in {\Gamma_1}$, and let $Q\subset\oG$ be an arbitrary Borel set; then the measures
\begin{equation*}
\begin{gathered}
\delta(y,Q)=\zeta(y)b_1(y)\chi_Q(\Omega_1(y)),\qquad
\alpha(y,Q)=b_2(y)\chi_Q(\Omega_2(y)),\\
\beta^1(y,Q)=\big(1-\zeta(y)\big)b_1(y)\chi_Q(\Omega_1(y))+\sum_{j=3,4}
b_j(y)\chi_Q(\Omega_j(y)),\qquad
 \beta^2(y,Q)=0
\end{gathered}
\end{equation*}
(for simplicity, we have omitted the subscript ``1'' in the notation of the measures)
satisfy Conditions~\ref{condK1}, \ref{cond1.2}, and \ref{cond2.3}--\ref{cond3.4}.

The author is grateful to Prof. A.L. Skubachevskii for attention
to this work.


\begin{thebibliography}{99}

\bibitem{Feller2}
{W. Feller},  {\it Trans. Amer. Math. Soc.}, {\bf 77}, 1--30
(1954).


\bibitem{Ventsel}
A. D. Ventsel,  {\it Teor. Veroyatnost. i Primen., \bf 4},
172--185 (1959); English transl.: {\it Theory Probab. Appl.}, {\bf
4} (1959).





\bibitem{SU}
K. Sato,  T.~Ueno,   {\it J. Math. Kyoto Univ.} {\bf 4}, 529--605
(1965).


\bibitem{BCP}
J.~M. Bony, P.~Courrege, P.~Priouret,   {\it Ann. Inst. Fourier
{\rm (}Grenoble{\rm )}} {\bf 18}, (1968) 369--521.

\bibitem{Taira3}
K. Taira, Boundary Value Problems, Feller Semigroups, and Markov
Processes, {\it in} ``Interaction between Functional Analysis,
Harmonic Analysis, and Probability'', Vol. 40, Marcel Dekker, New
York--Basel--Hong Kong, 1996.


\bibitem{Ishikawa}
Y. Ishikawa,   {\it J. Math. Soc. Japan}, {\bf 42}, 171--184
(1990).

\bibitem{GalSkJDE}
E. I. Galakhov, A. L. Skubachevskii,  {\it J.
      Differ. Equ.}, {\bf 176}, 315--355 (2001).

\bibitem{SkubDAN89}
A. L. Skubachevskii,   {\it Dokl. Akad. Nauk SSSR\/}, {\bf 307}, 287--291 (1989); English
     transl. in {\it Soviet Math. Dokl.\/}, {\bf 40} (1990).

\bibitem{GarMenaldi}
M. G. Garroni, J. L. Menaldi, {\it Second Order Elliptic
Integro-Differential Problems}, Chapman~$\&$~Hall/CRC, London--New
York--Washington, D.C., 2002.



\end{thebibliography}
\end{document}